\documentclass[a4paper,11pt]{article}
\usepackage{amssymb}
\usepackage{amsmath,amsfonts,amssymb,srcltx,amscd}
\usepackage[english]{babel}
\usepackage{indentfirst}
\usepackage{fancyhdr}
\usepackage{anysize}
\usepackage{xcolor}
\UseRawInputEncoding
\newtheorem{thm}{Theorem}[section]

\newcommand{\rev}[1]{\textcolor{blue}{#1}}

\begin{document}
\title{\rev{A Key Exchange Construction using Mihailova Subgroups in Braid groups}}

\author{Hanling Lin $^{1}$, Yu Han $^{1*}$ \\1.Shenzhen University, Shenzhen City 518060, China}
\date{}
\maketitle

\textbf{Abstract.} \rev{In this paper, we propose a modified Anshel-Anshel-Goldfeld (AAG) key exchange construction. The algebraic motivation underlying this construction comes from the membership problem for Mihailova subgroups of the braid group $B_n$ $(n \geq 6)$, a problem that is algorithmically unsolvable. We show that this perspective leads naturally to a quotient-group formulation involving Mihailova subgroups modulo the center of $B_n$. We also explain, however, that these algebraic facts do not by themselves provide a complete security proof for the protocol, because recovering a functionally equivalent conjugator modulo the center may already suffice for an adversary. Thus, the construction should be regarded as an algebraically motivated candidate whose full cryptographic security requires further study.}


\section{Introduction}\label{1}

With the rapid development of quantum computing, the security of many widely-used public-key cryptosystems, such as RSA and elliptic curve cryptography, is increasingly threatened. This urgent situation has motivated the cryptographic community to seek quantum-resistant alternatives, among which braid-based cryptography has emerged as a promising candidate. Early braid-based cryptographic protocols, notably those proposed by Anshel et al.\cite{AAG1} and Ko et al.\cite{Ko1}, rely primarily on the presumed difficulty of the conjugacy search problem (CSP) in braid groups. Unfortunately, this foundation has proven to be vulnerable: numerous attacks have been developed that successfully solve the CSP in braid groups, including length-based attacks (see \cite{GKTTV}, \cite{HS}, \cite{HT}, \cite{LL}, \cite{MU}), Burau representation attacks (see \cite{Hu}, \cite{LP}), super-summit set attacks \cite{FG}, and ultra-summit set attacks \cite{Ge}. As a consequence of these vulnerabilities, it has become widely accepted that the CSP alone cannot serve as a reliable security foundation for braid-based cryptographic schemes.

\rev{To address this fundamental security challenge, researchers have explored alternative mathematical problems within braid groups that might offer stronger security guarantees. A significant breakthrough came in 2006, when Shpilrain and Ushakov observed in \cite{SU2} that for $n \geq 6$, there exist certain subgroups (known as Mihailova subgroups) of the braid group $B_n$ for which the membership problem is algorithmically unsolvable---a fundamental result originally established by Mihailova in \cite{Mi}. Their discussion suggests that choosing private keys from such subgroups may impose an additional membership-type obstruction beyond the CSP. However, turning this observation into a complete cryptographic security proof requires further justification.}

Building upon this theoretical foundation, recent work has made the Mihailova subgroups more accessible for cryptographic applications. In \cite{WL}, Wang et al. provided an explicit presentation of the Mihailova subgroups of $F_2 \times F_2$, which possess an algorithmically unsolvable subgroup membership problem. Extending this result to braid groups, Lin et al.\cite{LWL} established, via an explicit isomorphism, that the braid group $B_n$ with $n \geq 6$ contains such Mihailova subgroups. These concrete characterizations provide the necessary tools to construct practical cryptosystems based on the membership problem.

\rev{In this paper, we use these recent developments to formulate a modified version of the Anshel-Anshel-Goldfeld (AAG) protocol from \cite{AAG1}, using the braid group $B_n$ with $n \geq 6$. Our key modification is to select private keys from Mihailova subgroups of $B_n$, while using the generators of $B_n$ as public keys. The main purpose of this work is to show that Mihailova subgroups can be incorporated naturally into an AAG-type construction and to analyze the algebraic consequences of doing so. In particular, we prove an unsolvability result for the corresponding quotient subgroup, while also explaining that this algebraic result alone does not imply resistance to all known attacks or a complete quantum-safe security proof.}

\rev{The remainder of this paper is organized as follows: In Section 2, we provide a comprehensive review of braid groups and their Mihailova subgroups, and present an explicit characterization of Mihailova subgroups of the braid group $B_n$ with $n \geq 6$ used in this paper. Section 3 reviews the original Anshel-Anshel-Goldfeld scheme and analyzes the known attacks~\cite{AAG1}. In Section 4, we present our modified key exchange scheme that incorporates Mihailova subgroups of the braid group. Finally, in Section 5, we discuss the algebraic implications and the present limitations of the associated security argument.}

\section{Braid groups and  Mihailova subgroups}\label{2}

For an integer $n \geq 2$, the braid group $B_{n}$ is defined by the following Artin presentation:
$$ B_{n}=\langle\sigma_{1},\sigma_{2},...,\sigma_{n-1}|\sigma_{i}\sigma_{j}\sigma_{i}=\sigma_{j}\sigma_{i}\sigma_{j}, \mbox{for} \, |i-j|=1,\sigma_{i}\sigma_{j}=\sigma_{j}\sigma_{i}, \mbox{for} \, |i-j|\geq2\rangle$$
where $\sigma_{1},\sigma_{2},...,\sigma_{n-1}$ are called the Artin generators of $B_{n}$ , and each element of $B_{n}$ is called an $n$-braid.

Let $H$ be a group defined by the presentation$$\Gamma=\langle x_{1}, x_{2}, ..., x_{k}|R_{1}, R_{2} ,..., R_{m}\rangle$$
where $k$ is an integer and $ k\geqslant2 $, and let $F_{k}$ be the free group on the generators $x_{1}, x_{2}, ..., x_{k}$. In the influential article \cite{Mi}, Mihailova associated to $H$ the Mihailova subgroup $M(H)$ of the direct product $F_{k}\times F_{k}$\textcolor{blue}{,} which is defined as
$$M(H)=\{(\omega_{1}, \omega_{2})|\omega_{1}=\omega_{2}, \omega_{1} \in H, \omega_{2} \in H\}$$
and proved the following theorem.

\begin{thm}\label{thm2.1} \emph{(Mihailova, \cite{Mi})} The membership problem for $M(H)$ in $F_{k}\times F_{k}$ is solvable if and only if the word problem for $H$ is solvable.
\end{thm}

In their study \cite{BV}, Bogopolski and Venura proved a theorem which gave an explicit representation of the Mihailova subgroup $M(H)$ in $F_{k}\times F_{k}$ if the group $H$ satisfies certain conditions. In \cite{WL}, Wang, Li and Lin gave a finite presentation of a group $H$, which is generated by only two elements and has an unsolvable word problem. Furthermore, they proved that the presentation of $H$ satisfies the conditions required in Bogopolski and Venura's theorem in \cite{BV}. Consequently, they gave an explicit presentation of the Mihailova subgroup $M(H)$ in $F_{2}\times F_{2}$ with a finite number of generators and infinitely countable defining relators as follows:
$$D=\langle(u, u), (t, t), (1, \mathcal{S}_{i}) |\mathcal{S}_{i}^{-1}(\delta^{-1}\mathcal{S}_{k}^{-1}\gamma_{k}^{-1}\delta)^{-1}\mathcal{S}_{i}(\delta^{-1}\mathcal{S}_{k}^{-1}\gamma_{k}^{-1}\delta), \,  \mathcal{S}_{i}^{-1}\gamma_{i}^{-1}\mathcal{S}_{i}\gamma_{i}\rangle$$
where $u$ and $t$ are the free generators of $F_2$,  $\delta  \in F_{2} \times F_{2}, \, \mathcal{S}_{i}=(\mathcal{R}_{i}^{(r)}(t,u))^{-1}\mathcal{R}_{i}^{(l)}(t, u), \, \mathcal{R}_{i}^{(r)}(t, u)$ and $\mathcal{R}_{i}^{(l)}(t, u)$ are defined in the Presentation C in \cite{WL}, $i$, $k = 1, 2, ..., 27$.

Since the word problem for $H$ is unsolvable, the main theorem in \cite{Mi} implies that the membership problem for the Mihailova subgroup $M(H)$ of $F_{2}\times F_{2}$ is also unsolvable.

For the braid group $B_{n}$ with $ n\geqslant6 $, Collins showed in \cite{Co} that the subgroups $G_i$ of $B_{n}$ generated by $\sigma_{i}^{2},\, \sigma_{i+1}^{2},\, \sigma_{i+3}^{2}$ and $\sigma_{i+4}^{2}$ (for $1\leq i \leq n-5$), denoted as $$G_{i}=\langle\sigma_{i}^{2},\sigma_{i+1}^{2},\sigma_{i+3}^{2},\sigma_{i+4}^{2}\rangle \quad(1\leq i \leq n-5),$$ are isomorphic to the direct product $F_{2}\times F_{2}$, where $F_{2}$ is the free group of rank 2.

It must also be mentioned that Lin et al.\cite{LWL} derived an explicit presentation for Mihailova subgroups of $B_n$ by introducing an isomorphism from $F_{2}\times F_{2}$ to the subgroups $G_{i}$ of $B_n$. Let $\phi$ be the isomorphism that maps the direct product $F_{2}\times F_{2}$, where $F_{2}$ is freely generated by $u$ and $t$, to the subgroup $G_{i}$ for $1\leq i \leq n-5$. This isomorphism $\phi$ is defined by
$$\phi:(u, 1)\mapsto\sigma_{i}^{2},\, (t, 1)\mapsto\sigma_{i+1}^{2}, \, (1, u)\mapsto\sigma_{i+3}^{2},\, (1, t)\mapsto\sigma_{i+4}^{2}$$

As a result, we can obtain a presentation of the Mihailova subgroup $M(G_{i})$ of $B_n$ and have the following theorem.

\begin{thm}\label{thm2.2} For the braid group $B_n$ with $ n\geqslant6 $, the subgroup membership problem for the Mihailova subgroups $M(G_{i})$(for $1\leq i \leq n-5$) of $B_{n}$ is unsolvable.
\end{thm}

\section{Review of Anshell-Anshell-Goldfeld scheme and Attacks}

\subsection{Anshell-Anshell-Goldfeld scheme}

The Anshel, Anshel, and Goldfeld key establishment protocol in \cite{AAG1} is recalled as follows. 

\noindent{\textbf{The public information:}} A group $G$\textcolor{blue}{,} two subgroups $S_A$ and $S_B$ of $G$ generated respectively by $\{s_1, \cdots, s_m\}$ and $\{t_1, \cdots, t_n\}$ with $s_1, \cdots, s_m, t_1, \cdots, t_n \in G$.

\noindent{\textbf{Key establishment protocol:}}
\begin{enumerate}
           \item [(1)] Alice selects an element $u=u(s_1, \cdots, s_m) \in S_A$ as her private key. Alice sends the elements $u^{-1}t_1 u, \, u^{-1}t_2 u, \cdots,\, u^{-1}t_n u$ to Bob.
           \item [(2)] Bob selects an element $v=v(t_1, \cdots, t_n) \in S_B$ as his private key. Bob sends the elements $v^{-1}s_1 v, \, v^{-1}s_2 v, \cdots, \,  v^{-1}s_m v$ to Alice.
           \item [(3)] Alice replaces $s_i$ with $v^{-1}s_{i}v$ in the word $u$ for $i=1, 2, \cdots, m$, and computes
               $$K_A=u^{-1}u(v^{-1}s_{1}v, v^{-1}s_{2}v, \cdots, v^{-1}s_mv)=u^{-1}v^{-1}uv$$
           \item [(4)] Bob replaces $t_{j}$ with $u^{-1}t_ju$ in the word $v$ for $j=1, 2, \cdots, n$, and computes
              $$K_B=(v^{-1}v(u^{-1}t_{1}u,\, u^{-1}t_{2}u, \, \cdots,\, u^{-1}t_nu))^{-1}=(v^{-1}u^{-1}vu)^{-1}=u^{-1}v^{-1}uv$$
\end{enumerate}

Clearly, $K_{A}=K_{B}$, which is therefore the shared key of Alice and Bob. The security is based on the difficulty of the conjugacy search problem in the base group $G$. As indicated in the papers (see \cite{HT}, \cite{LL}, \cite{HS}, \cite{GKTTV}), several algorithms have been found to address the conjugate search problem, which has led to significant challenges for cryptographic protocols based on such problems.

\subsection{Decision problems and Normal forms}

We list here the conjugacy search problem and the membership problem in groups, both of which serve as the theoretical foundation for numerous cryptosystems based on braid groups.

\textbf{Conjugacy Search Problem~(CSP)}:  Let $G$ be a group, and let $x, y$ be elements of $G$ satisfying $y=g^{-1}xg$ for some $g\in G$. The conjugate search problem is to find an element $g^{'}\in G$ such that $y=g^{'-1}xg^{'}$.

\textbf{Membership Problem~(MP)}: Let $G$ be a group, and let $H$ be a subgroup of $G$ generated by elements $b_{1}, b_{2}, ..., b_{k} \in G$. The subgroup membership problem for $H$ is, for any element $x$ of $G$, to determine whether $x \in H$, i.e., to determine whether $x$ can be expressed as a product of powers of $b_{1}, b_{2}, ..., b_{k}$.

Let $\Delta$ denote the fundamental braid $\Delta_{n}$ of the braid group $B_n$, which is defined inductively as follows:
$$  \Delta_1=1, \, \Delta_2=\Delta_1 \sigma_1, \, \cdots, \, \Delta_n=\Delta_{n-1}\sigma_{n-1}\sigma_{n-2}\cdots \sigma_2\sigma_1.$$
Then it is well known that for $n \geq 3$, the cyclic subgroup generated by $\Delta^2$ is the center of $B_n$ (see \cite{De} for reference).

Let $B_n^+$ denote the submonoid of $B_n$ generated by $\sigma_1, \cdots, \sigma_{n-1}$. Elements of $B_n^+$ together with the identity element $\varepsilon$ are called positive braids. A partial order $\leq$ on the elements of $B_n$ is defined by setting $u \leq v$ if and only if there exist positive braids $\alpha, \beta \in B_n^+$ such that $u=\alpha v \beta$. Any braid $\alpha \in B_n$ which satisfies $\varepsilon \leq \alpha \leq \Delta$ is called a canonical factor. 

Let $S_n$ be the symmetric group consisting of all permutations on the set $\{1, 2, \cdots, n\}$. There exists a canonical homomorphism $\pi : B_n \rightarrow S_n$ that maps each generator $\sigma_i$ to the transposition ($i$, $i + 1$), namely $\pi(\sigma_i)=(i, i+1)$. The restriction of $\pi$ to the set of canonical factors in $B_n$ induces a bijection \cite{em}.
A factorization $\gamma = \alpha \beta$ of a positive braid $\gamma$ into a canonical factor $\alpha$ and a positive braid $\beta$ is said to be left-weighted if and only if $\alpha$  attains the maximal length among all such decompositions. A right-weighted factorization is defined in an analogous manner.

We recall that every braid $w \in B_n$ can be written uniquely (see \cite{em,Gar}) as the following {\em  normal form}
\begin{equation}\label{Normalform}
   w = \Delta^k W_1 W_2 \cdots W_s,
\end{equation}
where each $W_i$ is a canonical factor and for all integers $1 \leq i < s$, the product $W_i W_{i+1}$ is left-weighted. In this normal form, the {\em infimum}, the {\em canonical length} and the {\em supremum} of $w$ are defined respectively as inf$(w)=k$, len$(w)=s$ and sup$(w)=k+s$.

 The explicit decomposition (\ref{Normalform}) is called the normal form of $w$, and $s$ is defined as the canonical length of $w$, denoted by $l(w)$. It is noted that each $W_i$ is a canonical factor, and thus can be uniquely represented by the corresponding permutation $\pi(W_i)$ in the symmetric group $S_n$. For a given braid word $w \in B_n$, its normal form can be computed in running time $O(|w|^2 n \log n)$ where $|\cdot|$ denotes the word length of $w$.

The super summit set $\mathcal{S}(w)$ of a braid $w$ is defined as the set of all conjugates of $w$ that attain the minimal possible canonical length. It can be proven that $\mathcal{S}(w)$ is finite; furthermore, there exists an algorithm to compute $\mathcal{S}(w)$ for any given braid word $w \in  B_n$ (see \cite{em}).

Let $w$ be a braid in $B_n$, and let $\Delta^k W_1 W_2 \cdots W_s$ be the normal form of $w$.  The braids $\partial_+(w)$ and $\partial_-(w)$ are defined as follows:
\begin{equation}\label{cycling}
  \partial_+(w)= \Delta^k W_2 \cdots W_s \tau^k(W_1), \, \,  \, \,   \partial_-(w)= \Delta^k \tau^k(W_s) W_1 W_2 \cdots W_{s-1}
\end{equation}
where $\tau$ is the flip automorphism that maps $\sigma_i$ to $\sigma_{n-i}$ for each $i$. We then say that
$\partial_+(w)$ (resp. $\partial_-(w)$) is obtained by cycling (resp. decycling) from $w$. Moreover, the braids $\partial_+(w)$ and $\partial_-(w)$ are both conjugates of $w$. So, one can construct the super summit set $\mathcal{S}(w)$ by means of the braids  $\partial_+(w)$ and $\partial_-(w)$.

\subsection{Attacks against Anshell-Anshell-Goldfeld scheme }

There are two primary kinds of attacks against the Anshell-Anshell-Goldfeld (AAG) scheme. One is length-based attacks \cite{GKTTV,HS,HT,LL,MU} and another is summit set-based attacks that include the super summit set attack \cite{em,Gar} and the ultra summit set attack \cite{Ge}. These attacks mainly aim at solving the conjugacy search problem (CSP).

Garside's approach to solving the Conjugacy Problem in the braid group $B_n$ mainly involves associating with every braid $b$ a distinguished finite set of its conjugates, called the Summit Set. El-Rifai and Morton showed in \cite{em} that the summit set can be replaced by one of its subsets, called the super summit set (SSS). Later, Gebhardt proposed a new refinement in \cite{Ge}, consisting in replacing the SSS with an even smaller set known as the ultra summit set (USS). Both the super summit set and ultra summit set are smaller and therefore easier to determine.

The common principle of these attacks is to try to retrieve a conjugator from a conjugate pair $(s, s^{'})$ by starting with $s^{'}$, which is assumed to be derived from $s$, and iteratively conjugating $s^{'}$ into a new braid $ts^{'}t^{-1}$ such that the length or the complexity of $ts^{'}t^{-1}$ is minimal. These attacks are particularly effective against the AAG scheme, where multiple pairs of conjugate braids associated with the same conjugating braid are known.

\rev{For protocols in which one of the public tuples is the standard generating tuple $(\sigma_1,\sigma_2,\cdots,\sigma_{n-1})$ of $B_n$, the multiple simultaneous conjugacy problem becomes especially relevant. Lee and Lee \cite{LL} proposed a summit-set based approach to such systems, Gonzalez-Meneses \cite{GM} improved algorithms for multiple simultaneous conjugacy problems in braid groups, and Myasnikov, Shpilrain, and Ushakov \cite{MSU} emphasized the cryptanalytic significance of highly structured public tuples. Since our protocol publishes the full tuple $(x^{-1}\sigma_1x,\cdots,x^{-1}\sigma_{n-1}x)$, these references must be taken into account in any security discussion of the present construction.}

\section{The reformed Anshell-Anshell-Goldfeld scheme}


\begin{thm}\label{thm1}
The subgroup $(G_i \langle\Delta^2\rangle)/\langle\Delta^2\rangle$ of the quotient group $B_n/\langle\Delta^2\rangle$ with $n \geq 6 $ and $1\leq i \leq n-5$ is isomorphic to $F_2 \times F_2$. Consequently, the subgroup 
$M(G_i\langle\Delta^2\rangle)/\langle\Delta^2\rangle$ is a Mihailova subgroup of $B_n/\langle\Delta^2\rangle$ and its membership problem is unsolvable.
\end{thm}

\noindent{Proof.}
Let $T_i=\langle \sigma_i^2, \sigma_{i+1}^2 \rangle$ be the subgroup of $G_i$ generated by $\{\sigma_i^2, \sigma_{i+1}^2 \}$. Then $T_i$ is isomorphic to the free group $F_2$ of rank 2. It follows that
$$(T_i\langle\Delta^2 \rangle)/\langle\Delta^2 \rangle = \langle \sigma_i^2\langle\Delta^2\rangle, \sigma_{i+1}^2\langle\Delta^2\rangle \rangle$$
which is the subgroup of $(G_i \langle\Delta^2\rangle)/\langle\Delta^2\rangle$ generated by $\{\sigma_i^2\langle\Delta^2\rangle, \sigma_{i+1}^2\langle\Delta^2\rangle \}$. To establish the desired result, it is sufficient to show that $T_i\cong (T_i\langle\Delta^2 \rangle)/\langle\Delta^2 \rangle$.  

Define a map $\theta :$ $\{\sigma_i^2, \sigma_{i+1}^2\}\mapsto\{\sigma_i^2\langle\Delta^2\rangle, \sigma_{i+1}^2\langle\Delta^2\rangle \}$ by $\theta(\sigma_i^2)=\sigma_i^2\langle\Delta^2\rangle, \theta(\sigma_{i+1}^2)=\sigma_{i+1}^2\langle\Delta^2\rangle$. Since $T_i$ is a free group of rank 2, this map induces a unique homomorphism from $T_i$ to $(T_i\langle\Delta^2 \rangle)/\langle\Delta^2 \rangle$. Let $g \in \ker(\theta)$. Then $g \in \langle\Delta^2\rangle$. Since $T_i$ is free and $\Delta^2$ commutes with every element of $T_i$, $g$ must be the identity element of $T_i$. Hence $\theta$ is an isomorphism.


The following key exchange scheme constitutes our reformation of the Anshel-Anshel-Goldfeld key exchange scheme proposed in~\cite{AAG1}.

\noindent{\textbf{The public information:}} 

$\cdot$ A braid group $B_n$ with $n \geq 6$ including the generators $\sigma_{1}$, $\sigma_{2}$, ... , $\sigma_{n-1}$;

$\cdot$ The Mihailova subgroups  $M_i=M(G_{i})$ with $1 \leq i \leq n-5$ as defined in Section~\ref{2}.

\noindent{\textbf{Key exchange phase:}}
\begin{enumerate}
\item [(1)] Alice selects a Mihailova subgroup $M_i=M(G_{i})$ for some integer $i$ ($1 \leq i \leq n-5$), and an element $x=x(\sigma_1, \sigma_2, \cdots, \sigma_{n-1})\in M_i$ in terms of the generators of $B_n$ as her private key. She then sends the elements $x^{-1}\sigma_{1}x, \, x^{-1}\sigma_{2}x, \, ...,\, x^{-1}\sigma_{n-1}x$ to Bob.
\item [(2)] Bob selects a Mihailova subgroup $M_j=M(G_{j})$ for some integer $j$ ($1 \leq j \leq n-5$) (where $M_j$ does not commute with $M_i$), and an element $y=y(\sigma_1, \sigma_2, \cdots, \sigma_{n-1}) \in M_j$ in terms of the generators of $B_n$ as his private key. He sends the elements $y^{-1}\sigma_{1}y, \, y^{-1}\sigma_{2}y, \,..., \,  y^{-1}\sigma_{n-1}y$ to Alice.
\item [(3)] Alice replaces $\sigma_{i}$ with $y^{-1}\sigma_{i}y$ in the word $x$ for $i=1, 2, \cdots, n-1$, and computes
               $$x(y^{-1}\sigma_{1}y, y^{-1}\sigma_{2}y, \cdots, y^{-1}\sigma_{n-1}y)=y^{-1}xy.$$
                She then multiplies $y^{-1}xy$ on the left by $x^{-1}$ to obtain the key $K_{A}=x^{-1}y^{-1}xy$;
\item [(4)] Bob replaces $\sigma_{j}$ with $x^{-1}\sigma_j x$ in the word $y$ for $j=1, 2, \cdots, n-1$, and computes
              $$K_{B}=(y^{-1}y(x^{-1}\sigma_{1}x,\, x^{-1}\sigma_{2}x, \, ...,\, x^{-1}\sigma_{n-1}x))^{-1}=x^{-1}y^{-1}xy$$

\end{enumerate}

Clearly, $K_{A}=K_{B}$; thus, the shared key for Alice and Bob is $K=K_{A}=K_{B}$.

\section{Parameters and Security Analysis}
{\color{blue}
The public tuple
\[
(x^{-1}\sigma_1x,\;x^{-1}\sigma_2x,\;\cdots,\;x^{-1}\sigma_{n-1}x)
\]
determines the corresponding inner automorphism of $B_n$. If an adversary finds an element $x' \in B_n$ such that
\[
x'^{-1}\sigma_i x' = x^{-1}\sigma_i x \qquad (i=1,2,\cdots,n-1),
\]
then $x'=C_xx$ for some $C_x \in B_n$, and hence
\[
x^{-1}C_x^{-1}\sigma_i C_xx=x^{-1}\sigma_i x \qquad (i=1,2,\cdots,n-1).
\]
Therefore $C_x$ commutes with every generator $\sigma_i$, so $C_x$ lies in the center $\langle \Delta^2 \rangle$ of $B_n$. The same argument applies to any element $y' \in B_n$ satisfying
\[
y'^{-1}\sigma_j y' = y^{-1}\sigma_j y \qquad (j=1,2,\cdots,n-1).
\]
Consequently, any such solutions have the form
\[
x'=\Delta^{2a}x,\qquad y'=\Delta^{2b}y
\]
for some integers $a,b$.

Since $\Delta^2$ is central, one has
\[
x'^{-1}y'^{-1}x'y'=x^{-1}y^{-1}xy.
\]
Thus, recovering a functionally equivalent conjugator modulo the center already suffices for an adversary to compute the same shared key. This shows that the original claim equating recovery of the public conjugators with membership in the quotient Mihailova subgroup was too strong.

Theorem~\ref{thm1} still yields a genuine algebraic consequence: the subgroup
\[
M(G_i\langle\Delta^2\rangle)/\langle\Delta^2\rangle
\]
of $B_n/\langle\Delta^2\rangle$ is a Mihailova subgroup with unsolvable membership problem. Hence, if an attacker were required to recover the exact cosets $x\langle\Delta^2\rangle$ and $y\langle\Delta^2\rangle$ together with subgroup-membership information, one would encounter an unsolvable membership obstruction. However, the protocol itself only requires the computation of a functionally equivalent conjugator modulo the center, and the gap between these two tasks is not bridged by the present argument.

Moreover, because the published conjugates are taken with respect to the full Artin generating tuple, summit-set and multiple simultaneous conjugacy methods such as those discussed in \cite{LL,GM,MSU} remain relevant to the cryptanalysis of the scheme. Consequently, the present construction should be regarded as an algebraically motivated candidate protocol rather than a fully proved secure key exchange scheme.
}

\section{Conclusion}

\rev{In this paper, we proposed an AAG-type key exchange construction in which the private keys are chosen from Mihailova subgroups of the braid group $B_n$. The main contribution is an algebraic analysis showing how unsolvability of the membership problem enters the quotient-group setting associated with the protocol. At the same time, the revised security discussion makes clear that this algebraic obstruction does not by itself constitute a complete proof of resistance to known cryptanalytic attacks. Establishing the full cryptographic security of the construction remains an open problem for future work.}

\end{document}